\documentclass[11pt]{article}
\usepackage{amssymb,amsmath,color} 

 \allowdisplaybreaks

\catcode`!=11
\let\!int\int \def\int{\displaystyle\!int}
\let\!lim\lim \def\lim{\displaystyle\!lim}
\let\!sum\sum \def\sum{\displaystyle\!sum}
\let\!sup\sup \def\sup{\displaystyle\!sup}
\let\!inf\inf \def\inf{\displaystyle\!inf}
\let\!cap\cap \def\cap{\displaystyle\!cap}
\let\!max\max \def\max{\displaystyle\!max}
\let\!min\min \def\min{\displaystyle\!min}
\let\!frac\frac \def\frac{\displaystyle\!frac}
\catcode`!=12

\let\oldsection\section
\renewcommand\section{\setcounter{equation}{0}\oldsection}

\allowdisplaybreaks

\newtheorem{thm}{Theorem}[section]
\newtheorem{pro}{Proposition}[section]

\newtheorem{re}{Remark}[section]

\begin{document}
\title{Global-in-time stability of 2D  MHD boundary layer in the Prandtl-Hartmann regime}
\author{{\bf Feng Xie}\\[1mm]
\small School of Mathematical Sciences, and LSC-MOE,\\[1mm]
\small Shanghai Jiao Tong University,
Shanghai 200240, P.R.China\\[1mm]
{\bf Tong Yang}\\[1mm]
\small Department of Mathematics,
City University of Hong Kong,\\[1mm]
\small Tat Chee Avenue, Kowloon, Hong Kong;\\
\small and School of Mathematical Sciences\\[1mm]
\small Shanghai Jiao Tong University,
Shanghai 200240, P.R.China
}
\date{}
\maketitle

\noindent{\bf Abstract:} In this paper, we prove global existence of
 solutions with analytic regularity to the 2D MHD boundary layer equations in
 the mixed Prandtl and Hartmann regime derived by formal multi-scale expansion in \cite{GP}. The analysis shows that the combined effect of the magnetic diffusivity and transveral magnetic field on the boundary leads to a linear damping on the tangential velocity field  near the boundary.
 And  this  damping effect  yields the global
 in time analytic norm estimate in the tangential space variable on the
 perturbation of  the classical steady Hartmann profile.

\footnotetext[1]{{\it E-mail address:}
tzxief@sjtu.edu.cn (F. Xie), matyang@cityu.edu.hk (T. Yang)}

\vskip 2mm


\vskip 2mm


\section{Introduction}
The following mixed Prandtl and Hartmann boundary layer equations from the classical incompressible MHD system were derived in \cite{GP} for flat boundary in two space dimensions (2D) when the physical parameters such as Reynolds number, magnetic Reynolds number and the Hartmann number satisfy some constraints in the high Reynolds numbers limit:
\begin{align}
\label{1.1}
\left\{
\begin{array}{ll}
\partial_tu_1+u_1\partial_xu_1+u_2\partial_yu_1=\partial_yb_1+\partial^2_yu_1,\\
\partial_yu_1+\partial^2_yb_1=0,\\
\partial_xu_1+\partial_yu_2=0, \quad x\in\mathbb{R},\,\, y\in\mathbb{R}_+.
\end{array}
\right.
\end{align}
Here, $(u_1,u_2)$ denotes the velocity field of the boundary layer and $b_1$ is the corresponding
tangential  magnetic component.

For the classic Hartmann bounday layer system, there is a family of
steady solutions called Hartmann layer. It turns out that Hartmann layer is
also a solution to the above system. In this paper, we will study the global-in-time stability of the Hartmann layer in the analytic function space. For this, consider the system
$(\ref{1.1})$ with  initial data
\begin{align}
\label{ID}
u_1(t=0,x,y)=u_{10}(x,y),
\end{align}
and the no-slip boundary conditions
\begin{align}
\label{1.2}
\begin{array}{ll}
u_1|_{y=0}=0,\quad
u_2|_{y=0}=0.
\end{array}
\end{align}
And the far field is  taken as a uniform constant state. Consequently, the pressure term vanishes in equation of $(\ref{1.1})_1$. Denote
\begin{align}
\label{1.3}
\lim_{y\rightarrow+\infty}u_1=\bar{u},\qquad \lim_{y\rightarrow+\infty}b_1=\bar{b}.
\end{align}
Integrating the equation of $(\ref{1.1})_2$ over $y$ yields
\begin{align}
\label{1.4}
-u_1(t,x,y)+\bar{u}=\partial_yb_1.
\end{align}
Thus, the equations (\ref{1.1}) can be written as
\begin{align}
\label{1.5}
\left\{
\begin{array}{ll}
\partial_tu_1+u_1\partial_xu_1+u_2\partial_yu_1=-u_1(t,x,y)+\bar{u}+\partial^2_yu_1,\\
\partial_xu_1+\partial_yu_2=0.
\end{array}
\right.
\end{align}
Recall that the classical Hartmann boundary layer is given by
\begin{align}
\label{1.6}
u_1=(1-e^{-y})\bar{u},\qquad u_2=0,
\end{align}
which  is a  steady solution to (\ref{1.5}).
Without loss of generality, set $\bar{u}=1$ and denote the perturbation by
$(u,v)$:
\begin{align}
\label{1.7}
u_1=(1-e^{-y})+u,\qquad u_2=v.
\end{align}
Obviously, $(u,v)$ satisfies
\begin{align}
\label{1.8}
\left\{
\begin{array}{ll}
\partial_tu+(1-e^{-y}+u)\partial_xu+v\partial_y(-e^{-y}+u)=-u+\partial^2_yu,\\
\partial_xu+\partial_yv=0,
\end{array}
\right.
\end{align}
with initial data
\begin{align}
\label{UI}
u_0(x,y)=u_{10}(x,y)-(1-e^{-y}),
\end{align}
and boundary condition
\begin{align}
\label{UB}
u|_{y=0}=0,\quad
v|_{y=0}=0.
\end{align}

Before we state the main result, let us first introduce the following weighted analytic regularity function spaces.
For some $r>1$, denote an analytic weight $M_m$  by
\begin{align*}
M_m=\frac{(m+1)^r}{m!}.
\end{align*}
With a parameter $\tau>0$, set
\begin{align}
&X_{m}=\|e^{\alpha y}\partial_x^m g\|_{L^2(\mathbb{R}^2_+)}\tau^{m}M_m,\quad Z_{m}=\|e^{\alpha y}\partial_y\partial_x^m g\|_{L^2(\mathbb{R}^2_+)}\tau^{m}M_m,\\
&Y_{m}=\|e^{\alpha y}\partial_x^m g\|_{L^2(\mathbb{R}^2_+)}\tau^{m-1/2}m^{1/2}M_m,\quad D_{m}=\|e^{\alpha y}\partial_x^m g\|_{L^\infty_yL^2_x}\tau^{m}M_m,
\end{align}
and define
\begin{align}\label{norm-1}
&\|g\|^2_{X^{r}_{\tau,\alpha}}=\sum_{m\geq 0}X^2_{m},\quad \|g\|^2_{Z^{r}_{\tau,\alpha}}=\sum_{m\geq 0}Z^2_{m},\\
&\|g\|^2_{Y^{r}_{\tau,\alpha}}=\sum_{m\geq 0}Y^2_{m},\quad\|g\|^2_{D^{r}_{\tau,\alpha}}=\sum_{m\geq 0}D^2_{m}.
\end{align}
Here, $\tau$ denotes the analytic radius.
\begin{thm}
\label{THM}
 Let the initial data $u_{10}(x,y)$ be a small perturbation of the Hartmann profile $(\bar{u}(1-e^{-y}),0)$ satisfying the compatibility conditions and
\begin{align}
\label{THM1}
\|\partial_yu_{10}+u_{10}-\bar{u}\|_{X^{r}_{\tau_0,\alpha}}\leq \delta_0,
\end{align}
with $r>1$, $0<\alpha<\sqrt{2}/2$ for
some small  constant $\delta_0>0$. Then there exists a unique global-in-time solution $(u_1, u_2)$ to the problem (\ref{1.1})-(\ref{1.3}) satisfying
\begin{align*}
\|g\|_{X^{r}_{\tau(t), \alpha}}\leq e^{-2(1-2\alpha^2)t}\delta_0,\quad
\mbox{with} \quad \tau(t) >\tau_0/2,
\end{align*}
for all time $t\ge 0$, where $g=\partial_yu_1+u_1-\bar{u}$.
\end{thm}
\begin{re}
The initial analytic radius $\tau_0$ and the size of initial perturbation $\delta_0$ should satisfy the constraint of (\ref{OT}).
\end{re}

Now let us briefly review the background and some related works. First of all,
the Prandtl equations derived by Prandtl \cite{P} in 1904 describe the
fluid phenomena near a boundary with no-slip boundary condition through the
high Reynolds number limit of the
incompressible Navier-Stokes equations. For this system, so far the mathematical theories are basically limited to two space dimensions except in the framework of
analytic or Gevrey function spaces or under some structure constraints. And
the justification of the Prandtl ansatz has been extensively investigated with
only a few results, cf. \cite{GN2,Mae,SC,SC1} and references therein. In fact,
under the monotonicity condition on the tangential velocity component
 in the normal direction, Oleinik firstly
obtained the local existence of classical solutions  in the two space
dimensions by using the  Crocco transformation, cf.  \cite{O}. This result
together with some other works in this direction are presented in
Oleinik-Samokhin's classical
book \cite{OS}. Recently, this well-posedness result was re-proved by using  energy method  in the framework of
weighted Sobolev spaces in \cite{AWXY} and \cite{MW1} independently, by observing the cancellation mechanism in the system.  By imposing an additional
favorable pressure condition,
a global in time weak solution was obtained by Xin and Zhang in
\cite{XZ}.

When this monotonicity condition is violated, seperation of the boundary
layer is  expected. For this, E-Engquist constructed a finite
time blowup solution to the Prandtl equations in \cite{EE}. And this kind of blowup result is extended to the van Dommelen-Shen type singularity  in \cite{KVW} when the outer Euler flow is spatially periodic. In addition, some interesting ill-posedness (or instability) phenomena of solutions to both
the linear and nonlinear Prandtl equations around
 a shear flow have been  studied, cf.
\cite{GD,GN,GE, GN1} and the survey paper \cite{E}.

 In the framework of analytic functions, Sammartino and Caflisch \cite{SC,SC1,CS} proved
the local well-posedness result of the Prandtl system and justified the
Prandtl ansatz. The analyticity requirement in the normal variable $y$ was later removed by
Lombardo, Cannone and Sammartino in \cite{LCS}. The main argument used in \cite{LCS,SC} is to apply
the abstract Cauchy-Kowalewskaya  theorem. For recent development
of mathematical theories in Gevrey function space to Prandtl equations, one can refer to \cite{GM,LY}.

A natural question then is  whether  global existence of smooth (or strong) solution can be achieved for the  Prandtl equations in either analytic or Gevrey regularity function spaces. However, the answer to this  is still not known. In this direction,
 a lower bound of  the life-span for small analytic solution to the  Prandtl equations with small perturbation analytic initial data was given by Zhang ang Zhang in \cite{ZZ}. Precisely, when  the outflow velocity is of the order
 of $\varepsilon^{5/3}$,  and the initial perturbation is  of the order of $\varepsilon$,  then the  Prandtl
system admits a unique analytic solution with life-span  greater than $\varepsilon^{-4/3}$. On the other hand, when the initial data is  a small perturbation
of a  Guassian error function,  almost global existence for the Prandtl  equations is obtained by Ignatova and Vicol in \cite{IV}, where the cancellation observed in \cite{MW1} and the monotonicity  of background solution are essentially used  to have a linear time decay damping effect.

Back to the MHD system,
it is believed that suitable magnetic field can stabilize the boundary layer in some physical regime \cite{R,A,D,G}. One can also refer to some very recent results in \cite{GP,LXY1,LXY2,LXY3} for the derivation of MHD boundary layer equations, stability analysis of magnetic field on the boundary layer from the mathematical point of view.

The purpose of this paper to show a global-in-time  existence of solution to a mixed Prandtl and Hartmann MHD boundary layer equations derived in \cite{GP}.
The key observation is that the combined effect of the magentic diffusivity and transveral magnetic field on the boundary leads to a linear damping on the tangential velocity field  near the boundary. This damping effect yields
 a time exponential decay in analytic regularity norm of the solution when
 we consider perturbation near the classical Hartmann profile.

Finally, the rest of the paper is organized as follows. We will reformulate the problem by using the cancellation mechanism observed in \cite{AWXY, MW1} for
the study of
Prandtl equations in Sobolev space in Section 2. The
uniform estimates on the solution in analytic norm will be given  in Section 3. Based on the uniform estimates, the global existence and uniqueness of solution to (\ref{1.1}) will be proved in the last section.
Throughout the paper,   $C$, $\bar{C}$, $C_0$ and $C_1$ are used to denote some
generic constants.

\section{Preliminaries}
As for the classical Prandtl equation, one needs to use the cancellation mechanism in the system to overcome the loss of derivative in order to
close the a priori estimate. For this, note that  the vorticity in the
boundary layer  $\omega=\partial_yu$ satisfies
\begin{align}
\label{1.9}
\partial_t\omega+(1-e^{-y}+u)\partial_x\omega-ve^{-y}+ v\partial_y\omega=-\omega+\partial^2_y\omega.
\end{align}
As in \cite{AWXY,MW1}, one can use the vorticity to cancel some term with
essential difficulty in the equation for $u$. Precisely,
by noticing the Hartmann layer $(u_1,0)$ has the
property $\frac{u_{1yy}}{u_{1y}}=-1$, set
\begin{align}
\label{1.10}
g=\omega+u,
\end{align}
then the new unknown function $g$ satisfies
\begin{align}
\label{1.11}
\partial_tg+(1-e^{-y}+u)\partial_xg+ v\partial_yg=-g+\partial^2_yg.
\end{align}
The relation between the new unknown function  $g$ and $u$ is
\begin{align}
\label{1.12}
u=e^{-y}\int_0^ye^zg(t,x,z)dz,
\end{align}
and  the initial data of $g$ is
\begin{align}
\label{IG}
g(0,x,y)=\partial_yu_{10}(x,y)+u_{10}(x,y)-1.
\end{align}
As for  the boundary condition on $g$, note that
 from (\ref{1.8}) and (\ref{UB}), we have
\begin{align*}
(\partial_y\omega-u)|_{y=0}=0,
\end{align*}
which implies
\begin{align}
\label{1.13}
(\partial_yg-g)|_{y=0}=0.
\end{align}
In the rest of the paper, we consider the reformulated
problem (\ref{1.11})-(\ref{1.13}).
 \begin{re}
\label{R2.1}
If  $\|g\|_{X^{r}_{\tau,\alpha}}<\infty$, then it follows from (\ref{1.12}) that $\|u\|_{X^{r}_{\tau,\alpha'}}<\infty$ and $\|u\|_{Z^{r}_{\tau,\alpha'}}<\infty$ with $0\leq \alpha'<\alpha$. In
addition, $\|u\|_{D^{r}_{\tau,\alpha}}<\infty$.
\end{re}
\section{Uniform Estimates}
In this section, we will derive the uniform estimates on the solution to (\ref{1.11})-(\ref{1.13}) in analytic regularity norm through energy  method.

Applying the operator $\partial_x^m$ on  (\ref{1.11}), multiplying the resulting equation by $e^{2\alpha y}\partial_x^mg$ and integrating it over $\mathbb{R}_+^2$
yield
\begin{align}
\label{2.1}
\int_{\mathbb{R}_+^2}\partial_x^m(\partial_tg+(1-e^{-y}+u)\partial_xg+ v\partial_yg+g-\partial^2_yg)e^{2\alpha y}\partial_x^mgdxdy=0.
\end{align}
We  estimate the above equation term by term. Firstly,
\begin{align}
\label{TE}
\int_{\mathbb{R}_+^2}\partial_x^m\partial_tge^{2\alpha y}\partial_x^mgdxdy=\frac12\frac{d}{dt}\|e^{\alpha y}\partial_x^mg\|^2_{L^2},
\end{align}
\begin{align}
\label{DE}
\int_{\mathbb{R}_+^2}\partial_x^mge^{2\alpha y}\partial_x^mgdxdy=\|e^{\alpha y}\partial_x^mg\|^2_{L^2},
\end{align}
and
\begin{align}
\label{2.2}
&-\int_{\mathbb{R}_+^2}\partial^2_y\partial_x^mge^{2\alpha y}\partial_x^mgdxdy\nonumber\\
=&\int_{\mathbb{R}}\partial_y\partial_x^mg(t,x,0)\partial_x^mg(t,x,0)dx+\|e^{\alpha y}\partial_y\partial_x^mg\|^2_{L^2}+2\alpha\int_{\mathbb{R}_+^2}\partial_y\partial_x^mge^{2\alpha y}\partial_x^mgdxdy\nonumber\\
=&\|\partial_x^mg(t,x,0)\|^2_{L^2_x}+\|e^{\alpha y}\partial_y\partial_x^mg\|^2_{L^2}-\alpha\|\partial_x^mg(t,x,0)\|^2_{L^2_x}-2\alpha^2\|e^{\alpha y}\partial_x^mg\|^2_{L^2}\nonumber\\
=&(1-\alpha)\|\partial_x^mg(t,x,0)\|^2_{L^2_x}+\|e^{\alpha y}\partial_y\partial_x^mg\|^2_{L^2}-2\alpha^2\|e^{\alpha y}\partial_x^mg\|^2_{L^2},
\end{align}
where in the second equality, we have used the boundary condition (\ref{1.13}).

For  the  two mixed nonlinear terms in (\ref{2.1}), we firstly have
\begin{align*}
&\int_{\mathbb{R}_+^2}\partial_x^m((1-e^{-y}+u)\partial_xg)e^{2\alpha y}\partial_x^mgdxdy\\
=&\sum_{j=0}^m(\begin{array}{ll}m\\j\end{array})\int_{\mathbb{R}_+^2}\partial^{m-j}_xu\partial_x^{j+1}ge^{2\alpha y}\partial_x^mgdxdy
\triangleq R_1,
\end{align*}
and
\begin{align*}
|R_1|\leq& \sum_{j=0}^{[m/2]}(\begin{array}{ll}m\\j\end{array})\|\partial^{m-j}_xu\|_{L^2_xL^\infty_y}\|e^{\alpha y}\partial^{j+1}_xg\|_{L^\infty_xL^2_y}\|e^{\alpha y}\partial^{m}_xg\|_{L^2}\\
+&\sum_{j=[m/2]+1}^{m}(\begin{array}{ll}m\\j\end{array})\|\partial^{m-j}_xu\|_{L^\infty_{xy}}\|e^{\alpha y}\partial^{j+1}_xg\|_{L^2}\|e^{\alpha y}\partial^{m}_xg\|_{L^2}.
\end{align*}
For $0\leq j\leq [m/2]$, by (\ref{1.12}), one has
\begin{align*}
&\|\partial^{m-j}_xu\|_{L^2_xL^\infty_y}=\|\partial^{m-j}_x\int_0^ye^{-(y-z)}g(t,x,z)dz\|_{L^2_xL^\infty_y}\\
\leq &\|e^{\alpha y}\partial^{m-j}_x\int_0^ye^{-(y-z)}g(t,x,z)e^{\alpha z}e^{-\alpha z}dz\|_{L^2_xL^\infty_y}\\
= &\|\int_0^ye^{-(1-\alpha)(y-z)}\partial^{m-j}_x g(t,x,z)e^{\alpha z}dz\|_{L^2_xL^\infty_y}
\leq C\|e^{\alpha y}\partial^{m-j}_x g\|_{L^2},
\end{align*}
provided that $\alpha<1$.
Using the Agmon inequality gives that
\begin{align*}
\|e^{\alpha y}\partial^{j+1}_xg\|_{L^\infty_xL^2_y}\leq C\|e^{\alpha y}\partial^{j+1}_xg\|^{1/2}_{L^2}\|e^{\alpha y}\partial^{j+2}_xg\|^{1/2}_{L^2}.
\end{align*}
For $[m/2]+1\leq j\leq m$,
\begin{align*}
&\|\partial^{m-j}_xu\|_{L^\infty_{xy}}
=\|\partial^{m-j}_x\int_0^ye^{-(y-z)}g(t,x,z)dz\|_{L^\infty_{xy}}\\
\leq &\|e^{\alpha y}\partial^{m-j}_x\int_0^ye^{-(y-z)}g(t,x,z)e^{\alpha z}e^{-\alpha z}dz\|_{L^\infty_{xy}}\\
= &\|\int_0^ye^{-(1-\alpha)(y-z)}\partial^{m-j}_x g(t,x,z)e^{\alpha z}dz\|_{L^\infty_{xy}}\\
\leq &C\|e^{\alpha y}\partial^{m-j}_x g\|_{L^2_yL^\infty_x}
\leq C\|e^{\alpha y}\partial^{m-j}_x g\|^{1/2}_{L^2}\|e^{\alpha y}\partial^{m-j+1}_x g\|^{1/2}_{L^2}.
\end{align*}
Consequently,
\begin{align}
\label{NE1}
&\sum_{m\geq 0}|R_1|\tau^{2m}M^2_m\nonumber\\
\leq& \frac{C}{(\tau(t))^{1/2}}\left\{\sum_{m\geq 0}\sum_{j=0}^{[m/2]}X_{m-j}Y_{j+1}^{1/2}Y_{j+2}^{1/2}Y_m
\frac{(\begin{array}{ll}m\\j\end{array})M_m}{M_{m-j}M_{j+1}^{1/2}M_{j+2}^{1/2}(j+1)^{1/4}(j+2)^{1/4}m^{1/2}} \right.\\
&\left.+\sum_{m\geq 0}\sum_{j=[m/2]+1}^{m}X_{m-j}^{1/2}X_{m-j+1}^{1/2}Y_{j+1}Y_m
\frac{(\begin{array}{ll}m\\j\end{array})M_m}{M_{m-j}^{1/2}M_{m-j+1}^{1/2}M_{j+1}(j+1)^{1/2}m^{1/2}}\right\}.\nonumber
\end{align}
The second nonlinear term can be estimated as follows. Note that
\begin{align*}
&\int_{\mathbb{R}_+^2}\partial_x^m(v\partial_yg)e^{2\alpha y}\partial_x^mgdxdy\\
=&\sum_{j=0}^m(\begin{array}{ll}m\\j\end{array})\int_{\mathbb{R}_+^2}\partial^{m-j}_xv\partial_x^{j}\partial_yge^{2\alpha y}\partial_x^mgdxdy
\triangleq R_2,
\end{align*}
and
\begin{align*}
|R_2|\leq& \sum_{j=0}^{[m/2]}(\begin{array}{ll}m\\j\end{array})\|\partial^{m-j}_xv\|_{L^2_xL^\infty_y}\|e^{\alpha y}\partial^{j}_x\partial_yg\|_{L^\infty_xL^2_y}
\|e^{\alpha y}\partial^{m}_xg\|_{L^2}\\
+&\sum_{j=[m/2]+1}^{m}(\begin{array}{ll}m\\j\end{array})\|\partial^{m-j}_xv\|_{L^\infty_{xy}}\|e^{\alpha y}\partial^{j}_x\partial_yg\|_{L^2}\|e^{\alpha y}\partial^{m}_xg\|_{L^2}.
\end{align*}
For $0\leq j\leq [m/2]$,
\begin{align*}
\|\partial^{m-j}_xv\|_{L^2_xL^\infty_y}=&\|\int_0^y(e^{-z}\int_0^ze^{s}\partial^{m-j+1}_xg(t,x,s)ds)dz\|_{L^2_xL^\infty_y}\\
=&\|\int_0^ye^{-\alpha z}(\int_0^ze^{(1-\alpha)(s-z)}\partial^{m-j+1}_xg(t,x,s)e^{\alpha s}ds)dz\|_{L^2_xL^\infty_y}\\
\leq& C\|e^{\alpha y}\partial^{m-j+1}_xg\|_{L^2},
\end{align*}
where we have used the fact that $0<\alpha<1$. Note that
\begin{align*}
\|e^{\alpha y}\partial^{j}_x\partial_yg\|_{L^\infty_xL^2_y}\leq C\|e^{\alpha y}\partial^{j}_x\partial_yg\|^{1/2}_{L^2}\|e^{\alpha y}\partial^{j+1}_x\partial_yg\|^{1/2}_{L^2}.
\end{align*}
For $[m/2]+1\leq j\leq m$,
\begin{align*}
\|\partial^{m-j}_xv\|_{L^\infty}=&\|\int_0^y(e^{-z}\int_0^ze^{s}\partial^{m-j+1}_xg(t,x,s)ds)dz\|_{L^\infty}\\
=&\|\int_0^ye^{-\alpha z}(\int_0^ze^{(1-\alpha)(s-z)}\partial^{m-j+1}_xg(t,x,s)e^{\alpha s}ds)dz\|_{L^\infty}\\
\leq& C\|e^{\alpha y}\partial^{m-j+1}_xg\|^{1/2}_{L^2}\|e^{\alpha y}\partial^{m-j+2}_xg\|^{1/2}_{L^2}.
\end{align*}
Consequently,
\begin{align}
\label{NE2}
&\sum_{m\geq 0}|R_2|\tau^{2m}M^2_m\nonumber\\
\leq& \frac{C}{(\tau(t))^{1/2}}\left\{\sum_{m\geq 0}\sum_{j=0}^{[m/2]}Y_{m-j+1}Z_{j}^{1/2}Z_{j+1}^{1/2}Y_m
\frac{(\begin{array}{ll}m\\j\end{array})M_m}{M_{m-j+1}M_{j}^{1/2}M_{j+1}^{1/2}(m-j+1)^{1/2}m^{1/2}} \right.\\
&\left.+\sum_{m\geq 0}\sum_{j=[m/2]+1}^{m}Y_{m-j+1}^{1/2}Y_{m-j+2}^{1/2}Z_{j}Y_m
\frac{(\begin{array}{ll}m\\j\end{array})M_m}{M_{m-j+1}^{1/2}M_{m-j+2}^{1/2}M_{j}(m-j+1)^{1/4}(m-j+2)^{1/4}m^{1/2}}\right\}.\nonumber
\end{align}
Combining the above estimates and using the definitions of $X^r_{\tau,\alpha}, Z^r_{\tau,\alpha}$ and $Y^r_{\tau,\alpha}$ give
\begin{align}
\label{OD}
&\frac12\frac{d}{dt}\|g\|^2_{X^r_{\tau, \alpha}}-\dot{\tau}\|g\|^2_{Y^r_{\tau, \alpha}}+\|g\|^2_{Z^r_{\tau, \alpha}}+(1-2\alpha^2)\|g\|^2_{X^r_{\tau, \alpha}}+(1-\alpha)\|g(t,x,0)\|^2_{X^r_{\tau, \alpha}}\nonumber\\
\leq &\sum_{m\geq 0}|R_1|\tau^{2m}M^2_m+\sum_{m\geq 0}|R_2|\tau^{2m}M^2_m.
\end{align}
Since
\begin{align}
\label{CM1}
&
\frac{(\begin{array}{ll}m\\j\end{array})M_m}{M_{m-j}M_{j+1}^{1/2}M_{j+2}^{1/2}(j+1)^{1/4}(j+2)^{1/4}m^{1/2}} \nonumber\\
= &\frac{(m+1)^r(j+1)(j+2)^{1/2}}{(m-j+1)^r(j+2)^{r/2}(j+3)^{r/2}(j+1)^{1/4}(j+2)^{1/4}m^{1/2}}\nonumber\\
\leq & C(1+j)^{1/2-r}
\end{align}
for all $0\leq j\leq [m/2]$, and
\begin{align}
\label{CM2}
&
\frac{(\begin{array}{ll}m\\j\end{array})M_m}{M_{m-j}^{1/2}M_{m-j+1}^{1/2}M_{j+1}(j+1)^{1/2}m^{1/2}}\nonumber\\
= &\frac{(m+1)^r(m-j+1)^{1/2}(j+1)^{1/2}}{(m-j+1)^{r/2}(m-j+2)^{r/2}(j+2)^{r}m^{1/2}}\nonumber\\
\leq & C(m-j+1)^{1/2-r},
\end{align}
for all $[m/2]+1 \leq j\leq m$.
Consequently,
\begin{align}
\label{E1}
\sum_{m\geq 0}|R_1|\tau^{2m}M^2_m\leq& \frac{C}{(\tau(t))^{1/2}}\left\{\sum_{m\geq 0}\sum_{j=0}^{[m/2]}X_{m-j}Y_{j+1}^{1/2}Y_{j+2}^{1/2}Y_m
(1+j)^{1/2-r}\right.\nonumber\\
&\left.+\sum_{m\geq 0}\sum_{j=[m/2]+1}^{m}X_{m-j}^{1/2}X_{m-j+1}^{1/2}Y_{j+1}Y_m
(m-j+1)^{1/2-r}\right\}\nonumber\\
\leq& \frac{C}{(\tau(t))^{1/2}}\|g\|_{X^{r}_{\tau, \alpha}}\|g\|^2_{Y^{r}_{\tau, \alpha}}.
\end{align}
Note that in the above second inequality,  we have used the following discrete Young's inequality
\begin{align*}
\|\zeta\cdot(\eta*\xi)\|_{l^1}\leq C\|\zeta\|_{l^2}\|\eta\|_{l^1}\|\xi\|_{l^2}
\end{align*}
with $\zeta_k=Y_k, \eta_k=Y^{1/2}_kY^{1/2}_{k+1}k^{1/2-r}, \xi_k=X_{k}$ for the first term on the right hand side
 in (\ref{E1}), and then H\"older inequality for
$
\|\eta\|_{l^1}\leq C\|g\|_{Y^{r}_{\tau, \alpha}}$,
provided
 $r>1$. And
for the second term on the right hand side of (\ref{E1}), we can choose
of $ \zeta_k=Y_k, \eta_k=X^{1/2}_kX^{1/2}_{k+1}(k+1)^{1/2-r}, \xi_k=Y_{k+1}$
with $
\|\eta\|_{l^1}\leq C\|g\|_{X^{r}_{\tau, \alpha}}$,
provided
 $r>1$. In conclusion, (\ref{E1}) holds if
$
r>1.
$

Similarly,
\begin{align}
\label{E2}
\sum_{m\geq 0}|R_2|\tau^{2m}M^2_m\leq& \frac{C}{(\tau(t))^{1/2}}\left\{\sum_{m\geq 0}\sum_{j=0}^{[m/2]}Y_{m-j+1}Z_{j}^{1/2}Z_{j+1}^{1/2}Y_m
(1+j)^{1/2-r}\right.\nonumber\\
&\left.+\sum_{m\geq 0}\sum_{j=[m/2]+1}^{m}Y_{m-j+1}^{1/2}Y_{m-j+2}^{1/2}Z_{j}Y_m
(m-j+1)^{1/2-r}\right\}\nonumber\\
\leq& \frac{C}{(\tau(t))^{1/2}}\|g\|_{Z^{r}_{\tau, \alpha}}\|g\|^2_{Y^{r}_{\tau, \alpha}},
\end{align}
provided that $r>1.$

It follows, from (\ref{OD}), (\ref{E1}) and (\ref{E2}), that
\begin{align}
&\frac12\frac{d}{dt}\|g\|^2_{X^{r}_{\tau, \alpha}}+\|g\|^2_{Z^{r}_{\tau, \alpha}}+(1-2\alpha^2)\|g\|^2_{X^{r}_{\tau, \alpha}}+(1-\alpha)\|g(\cdot,0)\|^2_{X^{r}_{\tau, \alpha}}\nonumber\\
\leq &(\dot{\tau}+\frac{C}{(\tau(t))^{1/2}}(\|g\|_{X^{r}_{\tau, \alpha}}+\|g\|_{Z^{r}_{\tau, \alpha}}))\|g\|^2_{Y^{r}_{\tau, \alpha}}.
\end{align}
Set
\begin{align}
\label{ODE}
\dot{\tau}+\frac{C}{(\tau(t))^{1/2}}(\|g\|_{X^{r}_{\tau, \alpha}}+\|g\|_{Z^{r}_{\tau, \alpha}})=0.
\end{align}
Then we have
\begin{align}
\frac12\frac{d}{dt}\|g\|^2_{X^{r}_{\tau, \alpha}}+\|g\|^2_{Z^{r}_{\tau, \alpha}}+(1-2\alpha^2)\|g\|^2_{X^{r}_{\tau, \alpha}}+(1-\alpha)\|g(\cdot,0)\|^2_{X^{r}_{\tau, \alpha}}
\leq 0.
\end{align}
When $
0<\alpha <\sqrt{2}/2$, we have
\begin{align}
\frac12\frac{d}{dt}\|g\|^2_{X^{r}_{\tau, \alpha}}+\|g\|^2_{Z^{r}_{\tau, \alpha}}+(1-2\alpha^2)\|g\|^2_{X^{r}_{\tau, \alpha}}
\leq 0.
\end{align}

 It follows that
\begin{align}
e^{2(1-2\alpha^2)t}\frac{d}{dt}\|g\|^2_{X^{r}_{\tau, \alpha}}+2e^{2(1-2\alpha^2)t}\|g\|^2_{Z^{r}_{\tau, \alpha}}+2(1-2\alpha^2)e^{2(1-2\alpha^2)t}\|g\|^2_{X^{r}_{\tau, \alpha}},
\leq 0,
\end{align}
that implies
\begin{align}
e^{2(1-2\alpha^2)t}\|g\|^2_{X^{r}_{\tau, \alpha}}+\int_0^t2e^{2(1-2\alpha^2)s}\|g(s)\|^2_{Z^{r}_{\tau, \alpha}}ds\leq \|g(0)\|^2_{X^{r}_{\tau_0, \alpha}}.
\end{align}

From  (\ref{ODE}), one has
\begin{align}
\label{AR}
&\tau(t)^{3/2}-\tau_0^{3/2}
=-C\int_0^t(\|g(s)\|_{X^{r}_{\tau, \alpha}}+\|g(s)\|_{Z^{r}_{\tau, \alpha}}))ds\nonumber\\
\geq &-C\int_0^t\|g(0)\|_{X^{r}_{\tau_0, \alpha}}e^{-(1-2\alpha^2)s}ds-C\int_0^t\|g(s)\|_{Z^{r}_{\tau, \alpha}}e^{(1-2\alpha^2)s}e^{-(1-2\alpha^2)s}ds\nonumber\\
\geq &-C_1\|g(0)\|_{X^{r}_{\tau_0, \alpha}}.
\end{align}
Hence, if the initial perturbation data is suitably  small such that
\begin{align}
\label{OT}
\frac{\tau_0}{K}> C_1^{2/3}\|g(0)\|^{2/3}_{X^{r}_{\tau_0, \alpha}},
\end{align}
with $K=(2\sqrt{2})^{2/3}/(2\sqrt{2}-1)^{2/3}$. Then (\ref{AR}) implies that
$
\tau(t) >\tau_0/2
$
for all $t\geq 0$.
Consequently, we have
\begin{pro}
\label{P3.1}
Under the same conditions of Theorem \ref{THM}, suppose that  $g$ is a solution to (\ref{1.11})-(\ref{1.13}) with analytic regularity in the norm
$X^{r}_{\tau, \alpha}$, then
\begin{align}
\label{P1}
e^{2(1-2\alpha^2)t}\|g\|^2_{X^{r}_{\tau, \alpha}}+\int_0^t2e^{2(1-2\alpha^2)s}\|g(s)\|^2_{Z^{r}_{\tau, \alpha}}ds\leq \|g(0)\|^2_{X^{r}_{\tau_0, \alpha}}
\end{align}
with $\tau(t) >\tau_0/2$ for all $t\geq 0$.
\end{pro}
\section{The Proof of Theorem \ref{THM}}

 By the uniform estimates on the solution to (\ref{1.11})-(\ref{1.13}) given
  in Proposition \ref{P3.1}, and the local existence of solutions in analytic
  function space that can be obtained by using the argument used in \cite{IV,ZZ},  the global existence of solution to (\ref{1.11})-(\ref{1.13}) follows.  Then
  by the relation (\ref{1.12}),  the global existence of $u$ to the initial-boundary value problem (\ref{1.8})-(\ref{UB}) is proved. In addition, according to Remark \ref{R2.1} and Proposition \ref{P3.1}, it follows that $\|u\|_{X^{r}_{\tau,\alpha'}}+\|\partial_yu\|_{X^{r}_{\tau,\alpha'}}<\infty$ with $0\leq \alpha'<\alpha$ and $\tau>\tau_0/2$. As consequence, the proof of global existence part in Theorem \ref{THM} is completed.

We now turn to prove the uniqueness. For this, it suffices to show the uniqueness of solution to (\ref{1.11})-(\ref{1.13}).
Assume that there are two solutions $g_1$ and  $g_2$ to (\ref{1.11})-(\ref{1.13}) with the same initial data $g_0$ satisfying $\|g_0\|_{X^{r}_{\tau_0,\alpha}}\leq \delta_0$. Denote the radii of analytic regularity for $g_1$ and  $g_2$ by $\tau_1(t)$ and $\tau_2(t)$ respectively.
Define $\tau(t)$ by
\begin{align}
\label{ODE-2}
\dot{\tau}+\frac{C}{(\tau(t))^{1/2}}(\|g_1\|_{X^{r}_{\tau_1, \alpha}}+\|g_1\|_{Z^{r}_{\tau_1, \alpha}})=0
\end{align}
with initial data
$
\tau(0)=\frac{\tau_0}{4}.
$
By the estimates given in Section 3, the analyticity radius $\tau(t)$ satisfies
\begin{align}
\label{U3}
\frac{\tau_0}{8}\leq \tau(t)\leq \frac{\tau_0}{4}\leq \frac{\min\{\tau_1(t),\tau_2(t)\}}{2},\quad \mbox{for}\quad t\geq 0.
\end{align}

Note that $\bar{g}=g_1-g_2$ satisfies
\begin{align}
\label{U4}
\partial_t\bar{g}+(1-e^{-y}+u_1)\partial_x\bar{g}+(v_1-v_2)\partial_yg_1=-\bar{g}+\partial_y^2\bar{g}+R,
\end{align}
with
\begin{align}
\label{U6}
R=-(u_1-u_2)\partial_xg_2-v_2\partial_y\bar{g}.
\end{align}
The initial data and the boundary condition for $\bar{g}$ are
\begin{align}
\label{DEI}
\bar{g}(t=0,x,y)=0,\quad
(\partial_y\bar{g}-\bar{g})|_{y=0}=0.
\end{align}
Following the arguments used  in Section 3, we have
\begin{align}
\label{UUU}
&\frac12\frac{d}{dt}\|\bar{g}\|^2_{X^{r}_{\tau, \alpha}}+\|\bar{g}\|^2_{Z^{r}_{\tau, \alpha}}+(1-2\alpha^2)\|\bar{g}\|^2_{X^{r}_{\tau, \alpha}}+(1-\alpha)\|\bar{g}(\cdot,0)\|^2_{X^{r}_{\tau, \alpha}}\nonumber\\
\leq &(\dot{\tau}+\frac{C}{(\tau(t))^{1/2}}(\|g_1\|_{X^{r}_{\tau, \alpha}}+\|g_1\|_{Z^{r}_{\tau, \alpha}}))\|\bar{g}\|^2_{Y^{r}_{\tau, \alpha}}
+\frac{C_0}{(\tau(t))^{1/2}}\|g_2\|_{Y^{r}_{\tau, \alpha}}(\|\bar{g}\|^2_{X^{r}_{\tau, \alpha}}+\|\bar{g}\|^2_{Z^{r}_{\tau, \alpha}}).
\end{align}
From (\ref{ODE-2}), we have
\begin{align}
\dot{\tau}+\frac{C}{(\tau(t))^{1/2}}(\|g_1\|_{X^{r}_{\tau, \alpha}}+\|g_1\|_{Z^{r}_{\tau, \alpha}})\leq 0,
\end{align}
where we have used the facts that $\tau(t)\leq \tau_1(t)$ and the norms $X^{r}_{\tau, \alpha}$ and $Z^{r}_{\tau, \alpha}$  are increasing with respect to $\tau$.
Moreover,
\begin{align}
\label{U13}
\frac{C_0}{(\tau(t))^{1/2}}\|g_2\|_{Y^{r}_{\tau, \alpha}}\leq& \frac{C}{\tau}\|g_2\|_{X^{r}_{2\tau, \alpha}}\leq \frac{C}{\tau}\|g_2\|_{X^{r}_{\tau_2, \alpha}}\nonumber\\
\leq& \frac{C_1}{\tau_0}\|g(0)\|_{X^{r}_{\tau_0, \alpha}}e^{-2(1-2\alpha^2)t}
\leq \bar{C}\delta_0^{1/3},
\end{align}
where in the second inequality $2\tau\leq \tau_2$ is used, and in the third inequality (\ref{OT}) and (\ref{P1}) are used.
By choosing $\delta_0$ suitably small, it follows, from (\ref{UUU})-(\ref{U13}), that
\begin{align}
\frac{d}{dt}\|\bar{g}\|^2_{X^{r}_{\tau, \alpha}}+\eta \|\bar{g}\|^2_{X^{r}_{\tau, \alpha}}\leq 0,
\end{align}
with $\eta$ being a small positive constant. This implies uniqueness of solution to (\ref{1.11}) for all $t>0$.
Then, the proof of uniqueness part in Theorem \ref{THM} is
completed.
\\

\noindent{\bf Acknowledgement:} Feng Xie's research is supported by National Nature Science Foundation of
China 11571231, the China Scholarship
Council and Shanghai Jiao Tong University SMC(A). Tong Yang's research is supported by National Nature Science Foundation of
China 11771169.

\end{document}